\documentclass{article}

\usepackage{zibtitlepage}
\usepackage[a4paper,left=.20\paperwidth,right=0.20\paperwidth]{geometry}

\usepackage{graphicx}

\usepackage{amsmath}
\usepackage{textcomp} 

\usepackage[mode=buildnew]{standalone}
\usepackage{tikz}
\usetikzlibrary{circuits}
\usepackage{pgfgas} 

\usepackage{natbib}

\input{commands}

\begin{document}
\ZTPAuthor{\ZTPHasOrcid{Felix Hennings}{0000-0001-6742-1983}}
\ZTPTitle{Large-scale empirical study on the momentum equation's inertia term}
\ZTPNumber{21-08}
\ZTPMonth{March}
\ZTPYear{2021}

\title{Large-scale empirical study on the momentum equation's inertia term}
\author{\ZTPHasOrcid{Felix Hennings}{0000-0001-6742-1983}}
\hypersetup{pdftitle={Large-scale empirical study on the momentum equation's inertia term},
  pdfauthor={Felix Hennings}}

\zibtitlepage
\maketitle

\begin{abstract}
A common approach to reduce the Euler equations' complexity for the simulation and optimization of gas networks is to neglect small terms that contribute little to the overall equations.
An example is the inertia term of the momentum equation since it is said to be of negligible size under real-world operating conditions.
However, this justification has always only been based on experience or single sets of artificial data points.
This study closes this gap by presenting a large-scale empirical evaluation of the absolute and relative size of the inertia term when operating a real-world gas network.
Our data consists of three years of fine-granular state data of one of the largest gas networks in Europe, featuring over 6,000 pipes with a total length of over 10,000 km.
We found that there are only 120 events in which a subnetwork consisting of multiple pipes has an inertia term of high significance for more than three minutes.
On average, such an event occurs less often than once every ten days.
Therefore, we conclude that the inertia term is indeed negligible for real-world transient gas network control problems.
\end{abstract}


\section{Introduction}
Solving the problem of optimizing the transient control of a gas transport network is a challenging task.
One source of complexity lies in the physical behavior of the gas, which is described for a one-dimensional pipe by a set of non-linear partial differential equations called the Euler equations, see \cite{Osi1996} for a general introduction.
A common approach from the literature for reducing the complexity of a corresponding optimization model is to ignore those terms in the equations, whose magnitudes are small compared to the other terms for variable values representing typical operating conditions.
This is, for example, done for the momentum equation, which can be expressed using mass flow \mFlow and pressure \press variables as
\begin{equation}
  \frac{1}{\area}\frac{\partial \mFlow}{\partial t} + \frac{\partial \press}{\partial x} + \frac{\sGasConst\temp}{\area^2}\frac{\partial (\mFlow^2 \zFact(\press)/\press)}{\partial x} + \frac{\fric \sGasConst \temp}{2 \area^2 \diam}\frac{|\mFlow| \mFlow \zFact(\press)}{\press} + \frac{\gravAcc \slope}{\sGasConst \temp}\frac{\press}{\zFact(\press)} =0. \label{eq:flowChgAnalysis_momentum_equation_raw}
\end{equation}
The other quantities are the space $x$ along the pipe, the time $t$, the cross-sectional area \area, diameter \diam, and slope \slope of the pipe, the specific gas constant \sGasConst, the gas temperature \temp, which we assume to be constant, the compressibility factor \zFact depending on the pressure, the darcy friction coefficient \fric, and the gravitational acceleration \gravAcc.
One usual assumption in the literature on gas network optimization is that the inertia term $\frac{1}{\area}\frac{\partial \mFlow}{\partial t}$ is small under usual or realistic operating conditions compared to the other terms of the momentum equation, especially the one describing the pipe wall friction $\frac{\fric \sGasConst \temp}{2 \area^2 \diam}\frac{|\mFlow| \mFlow \zFact(\press)}{\press)}$, see, for example, \cite{WilHolBatHan1964, FinGol1979, Osi1996, EhrSte2005, HerMohSac2010, BroGasHer2011, DomHilLanTis2017}.
To support this argument, they either refer to experience or derive the size of the different terms for one or two sets of exemplary values and show that the absolute size of the inertia term is small, usually less than 1\,\% of the overall sum of terms or the friction term.
However, to the best of our knowledge, there was no thorough empirical study on real-world data to challenge the given statement.
With this article, we close this gap.
Our evaluation is based on a large set of state data collected during our ongoing research cooperation with our project partner Open Grid Europe (\cite{OGE}), one of the major gas network operators in Europe.
This data enables us to evaluate the influence of the inertia term in real-world situations by determining its value, and measuring its absolute size as well as comparing it to the pipe wall friction term.

\section{Setup}
In preparation for the actual analysis, we have a more in-depth look at the used data set, the needed formulas, as well as the threshold definiting the relevance of the examined terms.

\subsection{Data}
Our study's data is provided by our project partner OGE\footnote{A special thanks to Julian Steinmeyer from OGE, as well as my colleagues Carsten Dre\ss{}ke and Tom Walther, for managing the compilation of the data archive.} and constitute a history of their network over a consecutive period of 41 months.
The network is represented as a directed graph and consists of roughly 8,000 nodes and 9,000 arcs representing the single elements of the network.
Of these, there are 6,000 pipes of a total length of more than 10,000\,km.
In the data, we are given the changing network topology as well as the state of the network for every 3 minutes during the time period.
This accumulates to roughly 3.6 billion or 3,600,000,000 data points combining the pipes with the points in time.
This state data contains, for example, values for the pressure and gas composition on all nodes, the flow on all arcs, and the currently used modes of the active elements.
Since the data is used for billing of the gas transport customers, it was checked and corrected multiple times, resulting in an overall high quality.

The data was created by a simulation-like process.
It starts from the current state of the network and is given the future values for the states of the active elements, the inflow values at the boundary nodes, the gas decomposition for the entry nodes, and a set of measured pressure profiles throughout the network.
The simulation software then determines those pressure, flow, and gas decomposition values for each network element, which fit the input data best.
This process is no pure simulation but a mixture of simulation and optimization as the problem is overdetermined by the given input.
Hence, the final result may differ from the given input values for some of the elements.
Especially for the measured pressure values, the difference may be rather high with up to 2 bar.

Despite the above-mentioned quality of the data, we found an entry node having erroneous data for its gas inflow for roughly three weeks.
Hence, we excluded the data of all those pipes from our analysis, which have been influenced by this entry during this period of time, resulting in 204,000 not evaluated data points. 

\subsection{Formulas}
To determine the value of the different terms in the momentum equation \eqref{eq:flowChgAnalysis_momentum_equation_raw}, we need to discretize it.
Therefore, we use the implicit Euler scheme for the time and integrate the equation along the length of the pipeline.
The remaining spacial integrals are approximated by the midpoint rule since we are only given a single average flow value per pipe in the data.
The result is the following discretization for a pipe $(\ell, r)$ and two consecutive time steps $t_0$ and $t_1$ with a time difference of \timeDiff seconds 
\begin{align}
  \pressI{\ell,t_1} - \pressI{r,t_1} =&\, \alpha + \beta + \gamma \label{eq:flowChgAnalysis_momentum_equation}\\
  \text{with} \qquad \alpha :=&\, \frac{\len}{\area\timeDiff}\left(\mFlowI{t_1} - \mFlowI{t_0}\right) \tag{inertia term} \\
  \beta :=&\, \frac{\fric \sGasConst \temp \len}{2 \area^2 \diam} \frac{|\mFlow_{t_1}| \mFlow_{t_1} \zFact(\pressI{m,t_1})}{\pressI{m,t_1}} \tag{friction term}\\
  \gamma :=&\, \frac{\sGasConst \temp}{\area^2} (\mFlow^2_{r, t_1} \frac{\zFact(\pressI{r,t_1})}{\pressI{r,t_1}} - \mFlow^2_{l, t_1} \frac{\zFact(\pressI{l,t_1})}{\pressI{l,t_1}}) + \frac{\gravAcc \slope \len}{\sGasConst \temp} \frac{\pressI{t_1}}{\zFact(\pressI{t_1})}. \tag{remaining terms}
\end{align}
We can interpret the equation as a definition of the pressure difference for the future time step $t_1$ given as the sum of the flow difference term $\alpha$, the friction term $\beta$, and the remaining terms $\gamma$.

We use the formula of \cite{Pap1968}(see \cite{Sal2002} for an english reference) for the compressibility factor \zFact and the one of Chen \cite{Che1979} for the darcy friction coefficient \fric.
For the pressure \pressI{m} at the midpoint, we use the average value of the pressures given at both end nodes, i.e., $\pressI{m} := \frac{\pressI{\ell}+\pressI{r}}{2}$.
Since we do not use the remaining terms in our following analysis, we do not need replacement values for the left flow \mFlowI{\ell} and right flow \mFlowI{r}.

In the data, flow is given as normal volumetric flow \nvFlow, which we will use in our analysis as well.
For the transformation of normal volumetric flow into mass flow \mFlow used in Equation~\eqref{eq:flowChgAnalysis_momentum_equation}, one multiplies by the gas density \nDens under normal conditions of 0\,\textdegree{}C and 1.01325\,bar pressure, equal to the pressure unit of 1 standard atmosphere.
The value \nDens depends on the gas decomposition and is given in the data set.
The final formulas for $\alpha$ and $\beta$ read as
\begin{align}
  \alpha &:= \frac{\len\nDens}{\area\timeDiff}\left(\nvFlowI{t_1} - \nvFlowI{t_0}\right) \label{eq:flowChgAnalysis_alpha_definition} \\
  \beta &:= \frac{\fric \sGasConst \temp \len \nDens^2}{2\area^2\diam} \frac{|\nvFlowI{t_1}|\nvFlowI{t_1}\zFact(\pressI{m,t_1})}{\pressI{m,t_1}}
\end{align}

\subsection{Thresholds} \label{sec:flowChgAnalysis_thresholdDefinition}
To determine the influence of the inertia term $\alpha$ in the momentum equation \eqref{eq:flowChgAnalysis_momentum_equation} at realistic operating conditions, we need to define certain thresholds for the term to be considered relevant or negligible for practical considerations.
The network of our project partner can be divided into areas featuring the majority of active elements in the network and long transport pipelines connecting these areas.
The pipelines typically have a length of 50\,km to 200\,km and cause a friction-induced pressure drop of several bar when transporting significant amounts of flow.
When comparing the solution found by an optimization algorithm to reality, there will always be small discrepancies between the calculated and the observed values caused by simplifying assumptions, incorrect input data, numerical problems, or other disturbances.
However, these small errors in the calculation can be compensated by slightly adjusting active elements' control when implementing the control recommendations.
Hence, we will focus on the potential errors caused in the single pipelines connecting these active elements.

Guided by the practitioners at OGE, we came up with the following absolute thresholds for practically relevant sizes of the inertia term $\alpha$ in the momentum equation for long pipeline sections between active elements:
\begin{align*}
  |\alpha|&<0.1\,\text{bar} && \text{No relevance} \\
  0.1\,\text{bar}\leq|\alpha|&<0.5\,\text{bar} && \text{Small relevance} \\
  0.5\,\text{bar}\leq|\alpha|&&& \text{High relevance}
\end{align*}
In addition to term' absolute size, we also demand it to have a minimum size compared to the other terms in the momentum equation to be relevant.
From these, the friction term $\beta$ is most of the time the dominant one.
Therefore, we define $\alpha$ to be negligible if $\frac{|\alpha|}{|\beta|}<0.01$ holds, which is in line with previous statements given in the literature, see, for example, \cite{Osi1996} or \cite{EhrSte2005}.

\section{Analysis}
Our analysis evaluates the statement that the inertia term $\alpha$ in the momentum equation \eqref{eq:flowChgAnalysis_momentum_equation} is small enough to be negligible under normal operating conditions.
Therefore, we will now look at the given data set and search for points in time at which there are pipes with inertia terms exceeding the negligibility thresholds defined above.

\subsection{Minimum flow change}
Our data set spans 41 months, consisting of about 1250 days in 3-minute granularity, making a total of 600,000 time points. 
Given the roughly 6000 pipes existing in the network, the entire data set contains about 3.6 billion data points representing the combination of pipes and time points.
As a first step, we aim at reducing this amount by looking only at those data points for which the flow value on the corresponding pipe changes by a minimal amount $\Delta_{\nvFlow}^\text{min}$.
This minimal amount should represent those data points having the potential to lead to significant values of $\alpha$, i.e., we choose $\Delta_{\nvFlow}^\text{min}$ largest possible such that all data points with $|\nvFlowI{t_1} - \nvFlowI{t_0}|<\Delta_{\nvFlow}^\text{min}$ are irrelevant based on the absolute threshold value of 0.1\,bar.
Hence, we are interested in
\begin{align*}
  \Delta_{\nvFlow}^\text{min} &:= \arg\max |\nvFlowI{t_1} - \nvFlowI{t_0}| \\
  \text{s.t.} \quad 0.1\,\text{bar} &> |\alpha| = \frac{\len\nDens}{\area\timeDiff}|\nvFlowI{t_1} - \nvFlowI{t_0}| \qquad \forall \len, \nDens, \area, \timeDiff \\
  \Leftrightarrow 0.1\,\text{bar} \frac{\area\timeDiff}{\len\nDens} &> |\nvFlowI{t_1} - \nvFlowI{t_0}| \qquad \forall \len, \nDens, \area, \timeDiff
\end{align*}
We deduce that $\Delta_{\nvFlow}^\text{min}$ is determined by $0.1\,\text{bar}\frac{\area^\text{min}\timeDiff^\text{min}}{\len^\text{max}\nDens^\text{max}}$ since
\begin{equation*}
  \frac{\area\timeDiff}{\len\nDens} > \frac{\area^\text{min}\timeDiff^\text{min}}{\len^\text{max}\nDens^\text{max}} \qquad \forall \len, \nDens, \area, \timeDiff.
\end{equation*}
For the corresponding extreme values, we choose the maximum length as $\len^\text{max} = 200\,\text{km}$ based on the above discussion, the minimal time interval as the one given in our data $\timeDiff^\text{min} = 3 \cdot 60\,\text{s}$, and a minimal $D^\text{min} = 150\,\text{mm}$ used for minimizing $A = D^2 \frac{\pi}{4}$.
There are, in fact, pipes in the network with a diameter of less than 150\,mm.
However, they are few and not used as transport pipelines connecting elements.
Furthermore, their accumulated length is less than 100\,km.
The typical natural gas mixtures occuring in the OGE network have a maximal normal density of $0.85\,\frac{\text{kg}}{\text{m}^3}$\cite{Cer2008}, which we round up to $\nDens^\text{max} = 0.9\,\frac{\text{kg}}{\text{m}^3}$.
Putting things together, we get
\begin{equation*}
  \Delta_{\nvFlow}^\text{min} = 0.1\,\text{bar}\frac{\area^\text{min}\timeDiff^\text{min}}{\len^\text{max}\nDens^\text{max}} = 0.636 \frac{1000\,\text{m}^3}{\text{h}}
\end{equation*}
To be on the safe side, we set $\Delta_{\nvFlow}^\text{min} = 0.5\,\frac{1000\,\text{m}^3}{\text{h}}$ for our analysis.
Using this value, we find about 730 million data points consisting of two subsequent time stamps and one pipe, such that the absolute flow difference on that pipe between the two timestamps exceeds $\Delta_{\nvFlow}^\text{min}$.

\subsection{Thresholds for single pipes}
In this second step, we actually calculate the values for $\alpha$ and $\beta$ for all 730 million data points determined in the previous step and evaluate both against the thresholds defined in Section~\ref{sec:flowChgAnalysis_thresholdDefinition}.
Note that $\alpha$'s absolute threshold is defined based on the assumption of a long pipeline with a length up to 200\,km.
Since all the single pipes in the network are much shorter, we determine for each pipe a normalized $\alpha$ value by dividing it by its length and comparing it against a length normalized threshold value.

\begin{figure}[ht]
  \centering
  \includegraphics[trim=0 10 0 50,clip,width=0.8\textwidth]{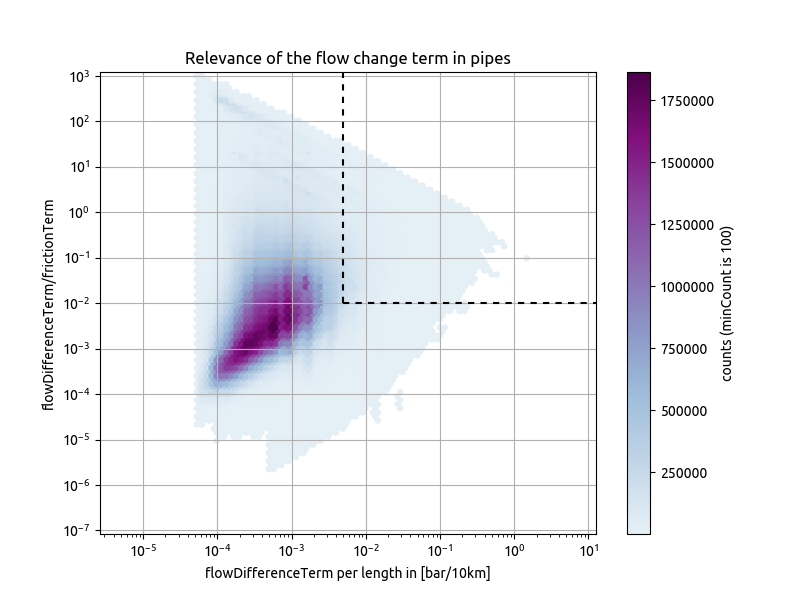}
  \caption{Hexagon plot of all data points having a flow change of at least $0.5\frac{1000\,\text{m}^3}{\text{h}}$. The color of the hexagon represents the number of associated original data points. The minimum count per hexagon is 100 ($\sim10^{-6}$ of all points). The dashed lines mark the area of values fulfilling the two threshold criteria for relevant points, i.e., all points outside of the marked top right area do not have relevant $\alpha$ values.}
  \label{fig:allFilteredRekoDataPoints}
\end{figure}
We plotted the results in Figure~\ref{fig:allFilteredRekoDataPoints}, using the normalized length $|\alpha|/10\,\text{km}$ and the friction value ratio $|\alpha|/|\beta|$ as axes.
The dashed lines mark the relevance thresholds, i.e., all points outside of the top right area do not have relevant $\alpha$ values.
For the length normalized $\alpha$ value, we will focus on points of small and high relevance having a minimum value of 0.1\,bar per 200\,km, which translates into a value of 0.005\,bar per 10\,km.

Figure~\ref{fig:allFilteredRekoDataPoints} shows that most points are clustered in the area between $\alpha$ values of $10^{-4}$ and $3\cdot 10^{-3}$ per 10\,km and $|\alpha|/|\beta|$ ratios of $10^{-4}$ and $10^{-1}$.
Due to the length normalized $\alpha$ threshold, all of these are not relevant for our ongoing analysis.
The diagonal line at the top right end of the point set indicates that extreme values regarding the two threshold criteria do not coincide.
Furthermore, we assume that the point set' straight-line boundary on the left side is caused by the minimal flow change threshold applied in the analysis' previous step.

In the following, we will continue to analyze those points fulfilling both relevance criteria, which are roughly 21 million data points. 

\subsection{Inertia term on paths}
In the previous steps, we have derived a set of pipe and time point combinations that potentially have relevant absolute sizes of $|\alpha|$ by normalizing their $\alpha$ values by the pipe length.
In this step, we will now look at the actual absolute $\alpha$ values produced by these data points.

First, we group the data points by time and look at each pair of subsequent time points individually.
For each of these pairs, we sort the data points into connected components based on the corresponding pipes and the underlying network topology.
In each connected component, we determine the $\alpha$-related absolute error as the length of a longest path in the component using the absolute $\alpha$ values as lengths per pipe.
The motivation for this procedure is the fact that the error caused by removing $\alpha$ from the momentum equation accumulates over space, along the pipes of the network.
In other words, the pressure difference of a path of pipes connecting two active elements is the sum of pressure differences defined by Equation~\eqref{eq:flowChgAnalysis_momentum_equation}, which includes the sum of corresponding $\alpha$ values for each of the pipes.

For creating the connected components, we group the pipes into sets connected directly by sharing an endnode or indirectly via other pipes in the set, valves being open at that point in time, or any resistors.
We do not use the active elements, i.e., regulators and compressor stations, to connect the pipes since we are interested in the error on pipelines connecting these elements.
As mentioned above, the maximum error inside a component is determined as the weight of a longest path when considering the absolute $\alpha$ values as edge weights.
Note that we use a directed graph as a representation of the connected component when searching for a longest path.
Each pipe is directed according to its topological orientation for positive $\alpha$ values and against the topological orientation otherwise. 
Open valves are added twice, once per orientation with a length value of zero.
For resistors, we are not guaranteed to have a flow value given in the data.
However, we can deduce the flow change from the pressure values of the endnodes: If the pressure reduction along the resistor increased, the flow through the resistor did so as well and we direct the arc according to its topological orientation.
If it decreased, we choose the opposite direction, and if it stayed the same, we add arcs for both directions.
In any of the cases, the weight of resistor arcs is zero.
An artificial example of a connected component with pipe flow values for a pair of time points, corresponding $\alpha$ values, and the longest path determining the $\alpha$ value for the whole component can be found in Figure~\ref{fig:flowChgAnalysis_alphaValuePaths}.

\begin{figure}[ht]
\centering

\newcommand{\nodeSquare}[2]{\draw[fill=white] (#1-0.05,#2-0.05) rectangle (#1+0.05,#2+0.05);}
\newcommand{\pynetPipe}[4]{\coordinate (LONGNAME_1) at (#1,#2);\coordinate (LONGNAME_2) at (#3,#4);\draw (LONGNAME_1) -- node[textNodeStyle,anchor=south] {huhu \makebox[1.0cm]{}} (LONGNAME_2);\draw[fill=white] ($(LONGNAME_1)!0.5!(LONGNAME_2)$) circle (0.08cm);}
\newcommand{\pynetPipeWithText}[6]{\coordinate (LONGNAME_1) at (#1,#2);\coordinate (LONGNAME_2) at (#3,#4);\draw (LONGNAME_1) -- node[textNodeStyle,anchor=#5] {#6} (LONGNAME_2);\draw[fill=white] ($(LONGNAME_1)!0.5!(LONGNAME_2)$) circle (0.08cm);}
\newcommand{\SB}{\rule[-0.1cm]{0.0cm}{0.1cm}} 
\begin{tikzpicture}
  [circuit,
   circuit symbol unit=2pt,
   text centered,
   textNodeStyle/.style={text width=0.5cm, text=black},
   arr/.style={->,thick,shorten <=3pt,shorten >=3pt,}]
  \scriptsize
  \draw (0,5) to [cs] (1,5) ;
  \pynetPipeWithText{1}{5}{2}{5}{south}{120 1200\SB};
  \draw (2,5) to [va] (3,5);
  \pynetPipeWithText{3}{5}{4}{ 6}{south east}{80 720};
  \pynetPipeWithText{3}{5}{4}{ 4}{north east}{40 360};
  \pynetPipeWithText{4}{ 6}{5}{ 6}{south}{80 640\SB};
  \pynetPipeWithText{4}{ 4}{5}{ 4}{north}{\makebox[1.0cm]{} 40 320};
  \pynetPipeWithText{5}{ 6}{6}{ 6}{south}{80 560\SB};
  \pynetPipeWithText{5}{ 4}{6}{ 4}{north}{\makebox[1.0cm]{} 40 280};
  \pynetPipeWithText{6}{ 6}{7}{5}{south west}{40 240};
  \pynetPipeWithText{6}{ 4}{7}{5}{north west}{40 240};
  \pynetPipeWithText{6}{ 6}{7}{7}{south east}{40 240};
  \pynetPipeWithText{7}{7}{8}{7}{south}{40 200\SB};
  \draw (8,7) to [cv] (9,7);
  \draw (7,5) to [re] (8,5);
  \pynetPipeWithText{8}{5}{9}{5}{south}{80 400\SB};
  \draw (9,5) to [cv] (10, 6);
  \draw (9,5) to [cv] (10, 4);
  \nodeSquare{0}{5}
  \nodeSquare{1}{5}
  \nodeSquare{2}{5}
  \nodeSquare{3}{5}
  \nodeSquare{4}{6}
  \nodeSquare{4}{4}
  \nodeSquare{5}{6}
  \nodeSquare{5}{4}
  \nodeSquare{6}{6}
  \nodeSquare{6}{4}
  \nodeSquare{7}{5}
  \nodeSquare{7}{7}
  \nodeSquare{8}{5}
  \nodeSquare{8}{7}
  \nodeSquare{9}{5}
  \nodeSquare{9}{7}
  \nodeSquare{10}{6}
  \nodeSquare{10}{4}
  \node at (5,3) {a) Network with flow value changes};
  \draw[arr,red] (1,0) to node[textNodeStyle,anchor=south] {0.09\SB} (2,0);
  \draw[arr,red] (2,0) to [bend left] node[textNodeStyle, anchor=south] {0.0\SB} (3,0);
  \draw[arr] (3,0) to [bend left] node[textNodeStyle, anchor=north] {\makebox[1.0cm]{} 0.0} (2,0);
  \draw[arr,red] (3, 0) to node[textNodeStyle, anchor=south east] {0.06} (4, 1);
  \draw[arr] (3, 0) to node[textNodeStyle, anchor=north east] {0.03} (4,-1);
  \draw[arr,red] (4, 1) to node[textNodeStyle, anchor=south] {0.05\SB} (5, 1);
  \draw[arr] (4,-1) to node[textNodeStyle, anchor=north] {\makebox[1.0cm]{} 0.03} (5,-1);
  \draw[arr,red] (5, 1) to node[textNodeStyle, anchor=south] {0.04\SB} (6, 1);
  \draw[arr] (5,-1) to node[textNodeStyle, anchor=north] {\makebox[1.0cm]{} 0.02} (6,-1);
  \draw[arr,red] (6, 1) to node[textNodeStyle, anchor=south west] {0.02} (7, 0);
  \draw[arr] (6,-1) to node[textNodeStyle, anchor=north west] {0.02} (7, 0);
  \draw[arr] (6, 1) to node[textNodeStyle, anchor=south east] {0.02} (7, 2);
  \draw[arr] (7, 2) to node[textNodeStyle, anchor=south] {0.01\SB} (8, 2);
  \draw[arr,red] (7, 0) to node[textNodeStyle, anchor=south] {0.0\SB} (8, 0);
  \draw[arr,red] (8, 0) to node[textNodeStyle, anchor=south] {0.03\SB} (9, 0);
  \nodeSquare{1}{0}
  \nodeSquare{2}{0}
  \nodeSquare{3}{0}
  \nodeSquare{4}{1}
  \nodeSquare{4}{-1}
  \nodeSquare{5}{1}
  \nodeSquare{5}{-1}
  \nodeSquare{6}{1}
  \nodeSquare{6}{-1}
  \nodeSquare{7}{0}
  \nodeSquare{7}{2}
  \nodeSquare{8}{0}
  \nodeSquare{8}{2}
  \nodeSquare{9}{0}
  \node at (5,-1.8) {b) Longest path in $\alpha$ value graph with total error value of 0.29};
\end{tikzpicture}
\caption{Example of finding a longest $\alpha$ value path in a connected component. The upper Figure a) shows a connected component enclosed by active elements, in which the numbers represent flow values for each pipe at the first point in time (upper value) and the subsequent time point (lower value). The lower Figure b) shows the directed arcs induced from the pipes, open valves, and resistors, where the numbers represent the $\alpha$ values on each of the arcs. The longest directed path regarding the $\alpha$ values represents the maximum error value in the component. Note that the pure sum of all $\alpha$ values in the component, as well as the undirected longest path in the component, would both overestimate the component's error value.
 }
\label{fig:flowChgAnalysis_alphaValuePaths}
\end{figure}

Note that the directed variant of the longest path search is necessary to prevent parallel loop lines in the component to count twice, i.e., once in both directions, even though their flow change and the associated $\alpha$ values are oriented towards the same node, see again Figure~\ref{fig:flowChgAnalysis_alphaValuePaths} for an example.
To calculate the longest path length, we search for a shortest path in the directed graph using negated weights for all pipes and negate the final result.
In the rare case of negative cycles, we save the negative cycle length, set the length of all affected arcs to zero, and add all negative cycle lengths to the final longest path length.
Hence, we may overestimate the length of the longest path and the connected component's $\alpha$ value.

When applying this procedure to the set of pipe and time point combinations, we find 6.3 million connected components created from the 21 million different data points.
From these connected components, 79,000 contain a longest path length in terms of absolute $\alpha$ values with length $\geq 0.1$\,bar and 2,500 have a path length of at least 0.5\, bar.
In other words, a flow change causing an absolute $\alpha$ value of more than 0.1\,bar on a path between active elements occurs on average every 23 minutes in the network.
A corresponding flow change for an absolute $\alpha$ value of at least 0.5\,bar happens on average every 12 hours.
In sum, 1.3 million single pipe data points are contained in the components for 0.1\,bar and 51,000 for 0.5\,bar, or 0.036\,\% and 0.0014\,\% of all 3.6 billion data points, respectively.
An overview of the number of components and data points for different thresholds can be found in Table ~\ref{tab:connectedComponentNumbers}.

\begin{table}[ht]
  \centering

  \begin{tabular}{rrrl} 
    Threshold & \#\,Components & \#\,Contained Pipes & avg comp occurance \\
    0.1 bar & 79,000 & 1,300,000 & 1 every 23 minutes\\
    0.2 bar & 19,000 & 370,000 & 1 every 1.6 hours\\
    0.3 bar & 7,500 & 150,000 & 1 every 4.0 hours \\
    0.4 bar & 3,900 & 81,000 & 1 every 7.7 hours\\
    0.5 bar & 2,500 & 51,000 & 1 every 12 hours\\
    0.6 bar & 1,700 & 35,000 & 1 every 18 hours\\
    0.7 bar & 1,200 & 26,000 & 1 every 25 hours\\
    0.8 bar & 930 & 20,000 & 1 every 1.3 days\\
    0.9 bar & 720 & 16,000 & 1 every 1.7 days\\
    1.0 bar & 580 & 13,000 & 1 every 2.2 days\\
  \end{tabular}
  \caption{Overview of connected component statistics for different thresholds of the absolute $\alpha$ value on a longest path in the component.
           For each threshold, we give the number of components having a longest path length of at least the given value and the sum of all pipe data points over all these components. In addition we give the average component occurance, determined from the number of components and the number of overall time steps. All Numbers are rounded to 2 significant digits.}
  \label{tab:connectedComponentNumbers}
\end{table}

\subsection{Inertia term over time} \label{sec:analyzeErrorsOverTime}
The previous section's results already indicate that $\alpha$ is indeed too small to be relevant for the vast majority of cases.
A connected component featuring highly relevant $\alpha$ values on average only occurs in the network every 12 hours, or each 240th time step, or 2500 times in total.
However, in this step of the analysis, we will have an even closer look at these remaining connected components and their 51,000 data points.

First, we determine how often in the 51,000 data points one pipe is contained in a highly relevant connected component at two or more subsequent time points.
If this is not the case, then the errors made by not including the $\alpha$ value in the pipe equation appear only for a brief period of time, i.e., that point in time in which the flow value changes to a different level.
More importantly, the $\alpha$ value errors decrease if we would increase the time granularity to more than 3 minutes.
In the literature on transient gas network optimization, the time step size usually ranges from 5 to 10 minutes, as in \cite{MakVanZloHij2016, BurEggGroMar2019}, to whole hours, like in \cite{Mor2007, DomGeiKolLan2011, ZloCheBac2015}.
If the $\alpha$ value is caused by a single large change in the flow values, its value would decrease by a factor of more than 3 for a 10 minute time granularity. 

In our data, 92\,\% of the 51,000 data points have not been in connected components with high $\alpha$ values for more than one consecutive time step.
Less than 6\,\% of the data points appear in connected components two times in a row, and less than 2\,\% of points in three components in a row.
The longest series of time points for a single pipe in a connected component with highly relevant $\alpha$ values was 8.
In total, there are only 8,700 out of the 51,000 data points for which the pipe is in a connected component of highly relevant $\alpha$ value also for the previous or next time interval.
These data points are contained in 570 different connected components in total.
Since we count a connected component per single time point, there are at most 285 consecutive time point series of length at least 2 having intersecting connected components of high $\alpha$ values.
Hence, these appear on average every 4.4 days.

\subsection{High flow change values}
As the final part of our analysis, we look at the size of absolute flow value changes in single data points.
From the formula \eqref{eq:flowChgAnalysis_alpha_definition} for $\alpha$ follows that it increases in absolute terms with an increasing absolute flow value difference $|\nvFlowI{t_1} - \nvFlowI{t_0}|$.
However, at a certain size, the flow changes are considered to be unrealistically high.
These values can occur in a simulation-like process like the one the data was generated with but are too large for real-world operating conditions.
Flow changes like these would put much strain on the single network elements, make the general network situation unstable, and cause a lot of noise, which all is not desirable from a network operator's point of view.
Therefore, there are elements to prohibit these large gas flow changes, for example, by slowly synchronize the pressure levels of two pipelines, which operated separately before and should now be connected.

For our analysis, we choose a rather high value of $\Delta_{\nvFlow}^\mathrm{real}=2,000,000\,\text{m}^3$/hour to be the limit for realistic flow value changes during a 3 minute time period.
This is roughly equal to the absolute size of the network's biggest entry being shutdown from maximum flow to zero during a 3 minutes time period.
Note that this value should actually be much smaller for pipes with smaller diameters since the corresponding velocity of the gas and hence also its noise level is higher for the same amount of flow.

From the 2,500 connected components of high relevance, roughly 800 contain a pipeline whose absolute flow change exceeds $\Delta_{\nvFlow}^\mathrm{real}$.
The remaining 1,700 connected components contain roughly 37,000 single data points.
If we repeated the time point series analysis from the previous Section~\ref{sec:analyzeErrorsOverTime}based on this reduced set of values, we find 4,400 data points in time point series of length at least two.
The longest time series are of length 4, i.e., present over 12 minutes, and all data points are contained in 240 connected components.
Therefore, there are at most 120 consecutive time point series of length at least 2 having intersecting connected components of high $\alpha$ values.
Hence, they appear on average every 10.4 days.
We consider this final average appearance rate of high relevant $\alpha$ values large enough to verify that $\alpha$ can be removed from the momentum equation since it is mostly irrelevant in realistic flow scenarios.

As a final remark, we want to highlight that after high-value flow changes, some parts of the change are often taken back during the following time steps.
This means that the flow stabilizes at some level, which is much closer to the original flow value than the one producing the large initial flow change value.
We did not quantify this effect, but when having unrealistically large flow changes, the following changes in the opposite direction can produce $\alpha$ high relevance values as well.
For example, we found out that every time series of length bigger than 5, produced from the set of 51,000 data points including the flow change value over $\Delta_{\nvFlow}^\mathrm{real}$, originates from two flow changes of 2,000,000\,m$^3$/hour and 10,000,000\,m$^3$/hour at one point in time and the following flow changes into the opposite direction during the subsequent time points.
These changes happening after the initial huge flow change would also be avoided in reality by reducing the size of the initial flow change and, therefore, further reducing the number of occurring network situations with relevant $\alpha$ values.
Having longer time steps would also reduce the $\alpha$ value in these situations since the initial flow changes are reduced by the following smaller changes in the opposite direction.

\section{Summary}

In this article, we evaluated the size of the inertia term $\alpha$ in the momentum equation on a large set of real-world gas network state data over a consecutive period of 41 months to review the statement that $\alpha$ is small and negligible under realistic operation conditions.
After establishing thresholds defining the relevant absolute sizes of the term and the relevant sizes of the term in relation to the friction value, we analyzed the data in multiple steps:
After reducing the number of data points to consider by sorting out those having a too-small flow change value, we found those data points for single pipes and time point pairs that have the potential to lead to relevant values of $\alpha$ when scaled to the maximum distance between active elements of 200\,km.
In the next step, we grouped these values by time and into connected components based on the network topology.
In these components, we established the overall absolute $\alpha$ value as the length of a longest directed path in terms of $\alpha$ values.
Components with highly relevant $\alpha$ values determined by this method occur on average every 12 hours in the network.
Finally, we determined those components in this set, which persist over more than a single 3-minute time step and do not feature unrealistically high flow values.
We found that there are only 120 of these components in the network during the entire time horizon, each lasting at most 12 minutes.
Hence they only appear less often than once every 10 days.
Therefore, we conclude that network situations featuring high relevant $\alpha$ values are rare enough to discard $\alpha$ from the momentum equation without losing much accuracy. 

\section*{Acknowledgements}
The work for this article has been conducted in the Research Campus MODAL funded by the German Federal Ministry of Education and Research (BMBF) (fund numbers 05M14ZAM \& 05M20ZBM).

%
%

\bibliographystyle{plainnat}
\bibliography{shortpaper}

\end{document}